\documentclass[12pt,leqno,twoside]{article}
\usepackage[latin1]{inputenc}
\usepackage[english]{babel}
\usepackage{amsthm,amssymb,amsmath,graphics,enumerate,epsfig,fancyhdr,url}

\newcommand{\ThisTitle}{}
\title{\ThisTitle}
\author{Magnus Dehli Vigeland\thanks{Department of Mathematics, University of Oslo, Norway.
{\it Email}\,: {\tt magnusv@math.uio.no}}
}
\date{}

\renewcommand{\ThisTitle}{Tropical complete intersection curves}

\theoremstyle{plain}
\newtheorem{theo}{Theorem}[section]

\newtheorem{lem}[theo]{Lemma}
\newtheorem{cor}[theo]{Corollary}

\theoremstyle{definition}
\newtheorem{dfn}[theo]{Definition}

\theoremstyle{remark}
\newtheorem{ex}[theo]{Example}
\newtheorem{remark}[theo]{Remark}

\newcommand{\pil}{\rightarrow}

\newcommand{\nn}{\mathbb{N}}
\newcommand{\zz}{\mathbb{Z}}
\newcommand{\pp}{\mathbb{P}}

\newcommand{\rr}{\mathbb{R}}
\newcommand{\La}{\Lambda}
\newcommand{\la}{\lambda}

\newcommand{\Ga}{\Gamma}
\newcommand{\ga}{\gamma}
\newcommand{\nsp}{\negmedspace}
\newcommand{\nnsp}{\negmedspace\negmedspace}

\newcommand{\cc}{\mathcal{C}}

\newcommand{\skriftt}{\mathcal{T}}
\newcommand{\skrifta}{\mathcal{A}}

\newcommand{\sub}{\subseteq}

\newcommand{\ijmn}{I_{\!J'}^{\,(0)}}
\newcommand{\ujn}{U_{\!J}^{\,(0)}}

\DeclareMathOperator{\conv}{conv}
\DeclareMathOperator{\Subdiv}{Subdiv}

\DeclareMathOperator{\vol}{vol}

\begin{document}

\maketitle
\begin{abstract}
A tropical complete intersection curve $\cc\sub \rr^{n+1}$ is a transversal intersection of $n$ smooth tropical hypersurfaces. We give a formula for the number of vertices of $\cc$ given by the degrees of the tropical hypersurfaces. We also compute the genus of $\cc$ (defined as the number of independent cycles of $\cc$) when $\cc$ is smooth and connected.
\end{abstract}

\section{Notation and definitions}
We work over the tropical semifield $\rr_{tr}\!=(\rr,\oplus,\odot)=(\rr,\max,+)$. A {\em tropical (Laurent) polynomial} in variables $x_1,\dotsc,x_m$ is an expression of the form
\begin{equation}\label{eq:f}
f=\nsp\bigoplus_{a=(a_1,\dotsc,a_m)\in\skrifta}\nsp\nsp\la_{a}\, x_1^{a_1}\dotsb x_m^{a_m}=\max_{a\in\skrifta}\{\la_{a}+ a_1x_1+\dotsb+a_mx_m\},
\end{equation}
where the {\em support set} $\skrifta$ is a finite subset of $\zz^m$, and the coefficients $\la_{a}$ are real numbers. (In the middle expression of \eqref{eq:f}, all products and powers are tropical.) The convex hull of $\skrifta$ in $\rr^m$ is called the {\em Newton polytope} of $f$, denoted $\Delta_f$.

Any tropical polynomial $f$ induces a regular lattice subdivision of $\Delta_f$ in the following way: With $f$ as in \eqref{eq:f}, let the {\em lifted Newton polytope} $\tilde{\Delta}_f$ be the polyhedron defined as
\begin{equation*}
\tilde{\Delta}_f:=\conv(\{(a,t)\:|\:a\in\skrifta, t\leq \la_{a}\})\;\sub \Delta_f\times \rr\sub\rr^m\times \rr
\end{equation*}
Furthermore, we define the {\em top complex} $\skriftt_f$ to be the complex whose maximal cells are the bounded facets of $\tilde{\Delta}_f$. Projecting the cells of $\skriftt_f$ to $\rr^m$ by deleting the last coordinate gives a collection of lattice polytopes contained in $\Delta_f$, forming a regular subdivision of $\Delta_f$. We denote this subdivision by $\Subdiv(f)$.

The standard volume form on $\rr^m$ is denoted by $\vol_m(\,\cdot\,)$, or simply $\vol(\,\cdot\,)$ if the space is clear from the context.

\subsection{Tropical hypersurfaces}\label{sec:tropbasics}
Note that any tropical polynomial $f(x_1,\dotsc,x_m)$ is a convex, piecewise linear function $f:\rr^m\pil\rr$.
\begin{dfn} Let $f:\rr^m\pil \rr$ be a tropical polynomial. The {\em tropical hypersurface} $V_{tr}(f)$ associated to $f$ is the non-linear locus of $f$.\end{dfn} 

It is well known that for any tropical polynomial $f$, $V_{tr}(f)$ is a finite connected polyhedral cell complex in $\rr^m$ of pure dimension $m-1$, some of whose cells are unbounded. Furthermore, $V_{tr}(f)$ is in a certain sense dual to $\Subdiv(f)$: There is a one-one correspondence between the $k$-cells of $V_{tr}(f)$ and the $(m-k)$-cells of $\Subdiv(f)$. A cell $C$ of $V_{tr}(f)$ is unbounded if and only if its dual $C^\vee\in\Subdiv(f)$ is contained in the boundary of $\Delta_f$. (For proofs consult \cite{Mikh_enum} and \cite{First}.)

Let $m\in\nn$, and let $e_1,\dotsc,e_m$ denote the standard basis of $\rr^m$. For any $d\in\nn_0$, we define the simplex $\Ga_d^m:=\conv\{0,de_1,\dotsc,de_m\}\sub\rr^m$, where $0$ denotes the origin of $\rr^m$.
For example, $\Ga_3^2$ is the triangle in $\rr^2$ with vertices $(0,0)$, $(3,0)$ and $(0,3)$. Note that
\begin{equation*}
  \vol(\Ga_d^m)=\frac1{m!}d^m.
\end{equation*}

\begin{dfn}
A tropical hypersurface $X=V_{tr}(f)\sub\rr^m$ is {\em smooth} if every maximal cell of $\Subdiv(f)$ is a simplex of volume $\frac1{m!}$. If in addition we have $\Delta_f=\Gamma_d^m$ for some $d\in\nn$, we say that $X$ is smooth of {\em degree} d.  
\end{dfn}

\subsection{Minkowski sums and mixed subdivisions}\label{sec:mink}
The set $\mathcal{K}^m$ of all convex sets in $\rr^m$ has a natural structure of a semiring, as follows: If $K_1$ and $K_2$ are convex sets, we define binary operators $\oplus$ and $\odot$ by
\begin{align}
K_1\oplus K_2&:=\conv(K_1\cup K_2)\\
\label{mink} K_1\odot K_2&:=K_1 + K_2.
\end{align}
The operator $+$ in \eqref{mink} is the {\em Minkowski sum}, defined for any two subsets $A,B\sub\rr^m$ by
$  A + B:=\{a+b\:|\:a\in A,b\in B\}$.
The Minkowski sum of two convex sets are again convex, so \eqref{mink} is well defined. Furthermore, it is easy to see that $\odot$ distributes over $\oplus$, and it follows that $\mathcal{K}^m$ is indeed a semiring.

\begin{lem}\label{lem:semihomo}
Let $\rr_{tr}[x_1,\dotsc,x_m]$ be the semiring of tropical polynomials in $n$ variables. The map $\rr_{tr}[x_1,\dotsc,x_m]\pil \mathcal{K}^{m+1}$ defined by $f\mapsto \tilde{\Delta}_f$, is a homomorphism of semirings.
\end{lem}
\begin{proof}
This is a straightforward exercise. The key ingredients are the identities 
\begin{equation*}\begin{split}
\conv(A\cup B)&=\conv(\conv(A)\cup \conv(B))\quad\text{and}\\
\conv(A + B)&=\conv(A) + \conv(B),
\end{split}\end{equation*}
which hold for any (not necessarily convex) subsets $A,B\sub \rr^m$.
\end{proof}

Let $f_1,\dotsc,f_n$ be tropical polynomials, and set $f:=f_1\odot\dotsb\odot f_n$. As a consequence of Lemma \ref{lem:semihomo}, we find that $\Subdiv(f)$ is the subdivision of $\Delta_f=\Delta_{f_1}+\dotsb +\Delta_{f_n}$ obtained by projecting the top complex of $\tilde{\Delta}_f=\tilde{\Delta}_{f_1}+ \dotsb +\tilde{\Delta}_{f_n}\sub \rr^m\times\rr$ to $\rr^m$ by deleting the last coordinate. 

For any cell $\La\in \Subdiv(f)$, the lifted cell $\tilde{\La}\in \skriftt_f$ can be written uniquely as a Minkowski sum $\tilde{\La}=\tilde{\La}_1 +\dotsb + \tilde{\La}_n$, where $\tilde{\La}_i\in \skriftt_{f_i}$ for each $i$. Projecting each term to $\rr^m$ gives a representation of $\La$ as a Minkowski sum $\La=\La_1 +\dotsb + \La_n$. The subdivision $\Subdiv(f)$, together with the associated Minkowski sum representation of each cell, is called the regular {\em mixed subdivision} of $\Delta_f$ induced by $f_1,\dotsc,f_n$. 

\begin{remark}
Note that the representation of $\La$ as a Minkowski sum of cells of the $\Subdiv(f_i)$'\,s is not unique in general. Following \cite{bert}, we call the representation obtained from the lifted Newton polytopes as described above, the {\em privileged representation} of $\La$.
\end{remark}

\begin{dfn}
The {\em mixed cells} of the mixed subdivision are the cells with privileged representation $\La=\La_1 +\dotsb + \La_n$, where $\dim \La_i\geq 1$ for all $i=1,\dotsc, n$. 
\end{dfn}
\section{Intersections of tropical hypersurfaces}
In this section we go through some basic properties and definitions regarding unions and intersections of tropical hypersurfaces. Most of the material here also appear in the recent article \cite{bert}.

We begin by observing that any union of tropical hypersurfaces is itself a tropical hypersurface. This follows by inductive use of the following lemma:
\begin{lem}\label{union}
If $X$ and $Y$ are tropical hypersurfaces in $\rr^m$, and $f,g$ are tropical polynomials such that $X=V_{tr}(f)$ and $Y=V_{tr}(g)$, then $X\cup Y=V_{tr}(f\odot g)$.
\end{lem}
\begin{proof}
By definition, $V_{tr}(f\odot g)$ is the non-linear locus of the function $f\odot g=f+g$. Since $f$ and $g$ are both convex and piecewise linear, this is exactly the union of the non-linear loci of $f$ and $g$ respectively.
\end{proof}

\begin{remark}\label{rem:subu}
Let $U=X_1\cup\dotsb\cup X_n$, where $X_i=V_{tr}(g_i)\sub \rr^m$ is a tropical hypersurface for each $i$. We denote by $\Subdiv_U$ the mixed subdivision of $\Delta_{g_1}+\dotsb +\Delta_{g_n}$ induced by $g_1,\dotsc, g_n$. It follows from Lemma \ref{union} and the discussion in Section \ref{sec:mink} that $\Subdiv_U$ is dual to $U$ in the sense explained in Section \ref{sec:tropbasics}.
\end{remark}

Moving on to intersections, we will only consider smooth hypersurfaces. Let $I$ be the intersection of smooth tropical hypersurfaces $X_1,\dotsc,X_n\sub\rr^m$, where $n\leq m$. As a first observation, notice that $I$ is a polyhedral complex, since the $X_i$'s are. The intersection is {\em proper} if $\dim(I)=m-n$. 

Let $C$ be a non-empty cell of $I$. 
Then $C$ can be written uniquely as $C=C_1\cap \dotsb \cap C_n$, where for each $i$, $C_i$ is a cell of $X_i$ containing $C$ in its relative interior. (The relative interior of a point must here be taken to be the point itself.) 

Regarding $C$ as a cell of the union $U=X_1\cup\dotsb\cup X_n$, we consider the dual cell $C^\vee\in\Subdiv_U$ (cf. Remark \ref{rem:subu}). From Section \ref{sec:mink}, we know that $C^\vee$ has a privileged representation as a Minkowski sum of cells of the subdivisions dual to the $X_i$'s. It is not hard to see that this representation is precisely $C^\vee=C^\vee_1+ \dotsb + C^\vee_n$. In particular, since $\dim C_i\leq m-1$, and therefore $\dim C^\vee_i\geq 1$, for each $i$, $C^\vee$ is a mixed cell of $\Subdiv_U$.

\begin{dfn}\label{transdef}
With the notation as above, the intersection $X_1\cap\dotsb \cap X_n$ is {\em transversal along $C$} if
\begin{equation}\label{tight}
\dim C^\vee=\dim C^\vee_1+\dotsb +\dim C^\vee_n.
\end{equation}
More generally, the intersection $X_1\cap \dotsb\cap X_n$ is said to be {\em transversal} if for any subset $J\sub\{1,\dotsc,n\}$ (of size at least two), the intersection $\bigcap_{i\in J} X_i$ is proper and transversal along each cell. 
\end{dfn}

\begin{remark}
Definition \ref{transdef} implies that if smooth tropical hypersurfaces $X_1,\dotsc,X_n$ intersect transversely, then $\Subdiv_U$ is a {\em tight coherent mixed subdivision} (see e.g. \cite{SturmViro}).
\end{remark}

Recall from standard theory that the {\em $k$-skeleton} $X^{(k)}$ of a polyhedral complex $X$, is the subcomplex consisting of all cells of dimension less or equal to $k$. It is not hard to see from Definition \ref{transdef} that if $X$ and $Y$ are tropical hypersurfaces intersecting transversely in $\rr^n$, then
\begin{equation}\label{xytrans}
X^{(j)}\cap Y^{(k)}=\emptyset
\end{equation}
for all nonnegative integers $j,k$ such that $j+k<n$. More generally, we find that:
\begin{lem}
Suppose $X_1,\dotsc, X_n$ intersect transversally, and let $I_J=\bigcap_{i\in J} X_i$, where $J$ is a subset of $\{1,2,\dotsc, n\}$. For each $s\notin J$ we have
\begin{equation*}
I_J^{(j)}\cap X_s^{(k)}=\emptyset,
\end{equation*}
for all $j,k$ such that $j+k<n$.
\end{lem}

\begin{figure}[tbp]
\begin{minipage}[b]{.48\linewidth}
\centering
\includegraphics[height=4cm]{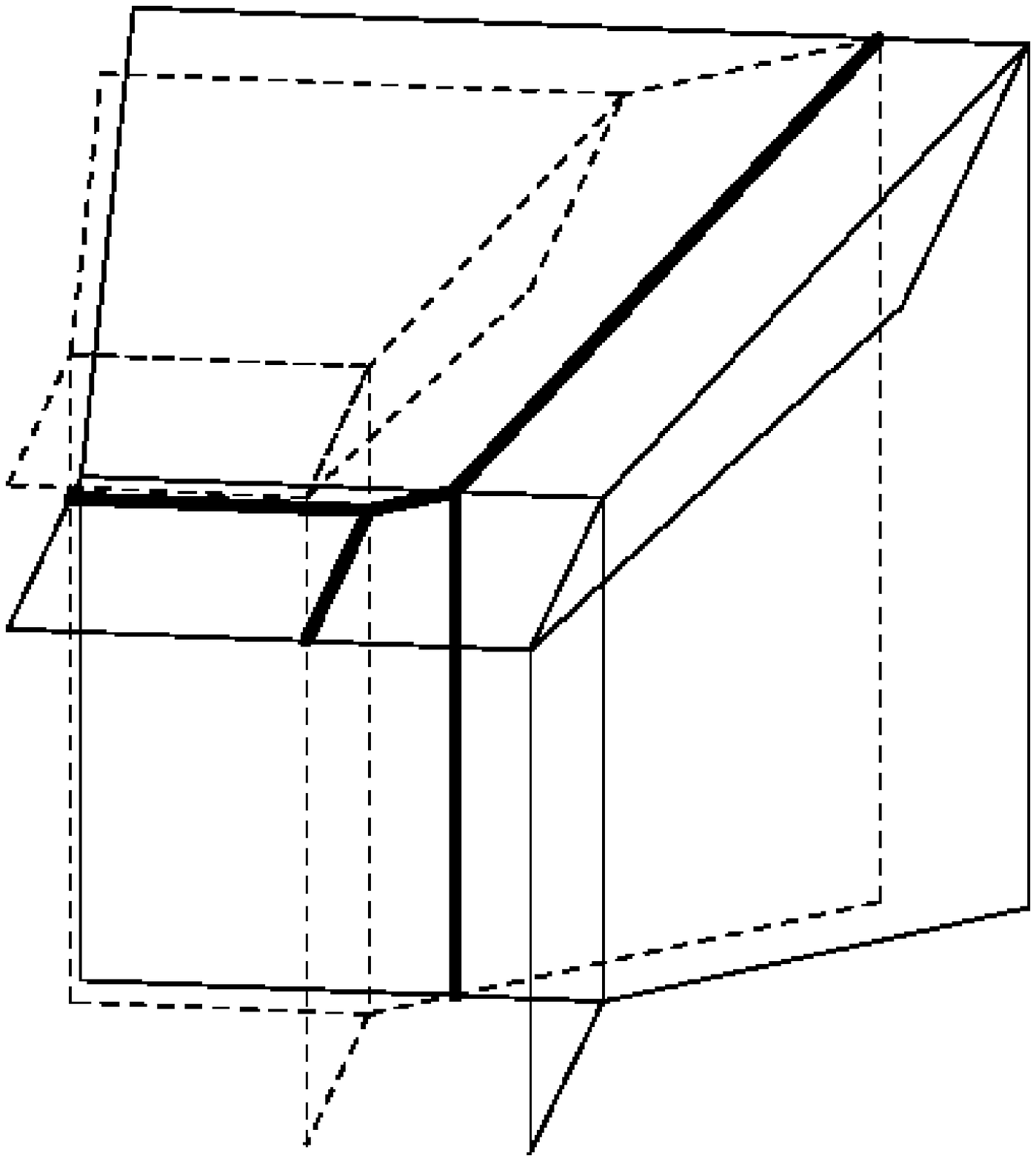}
\caption{Tropical planes intersecting in a tropical line.}\label{3line}
\end{minipage}\hfill
\begin{minipage}[b]{.48\linewidth}
  \centering
\includegraphics[height=4cm]{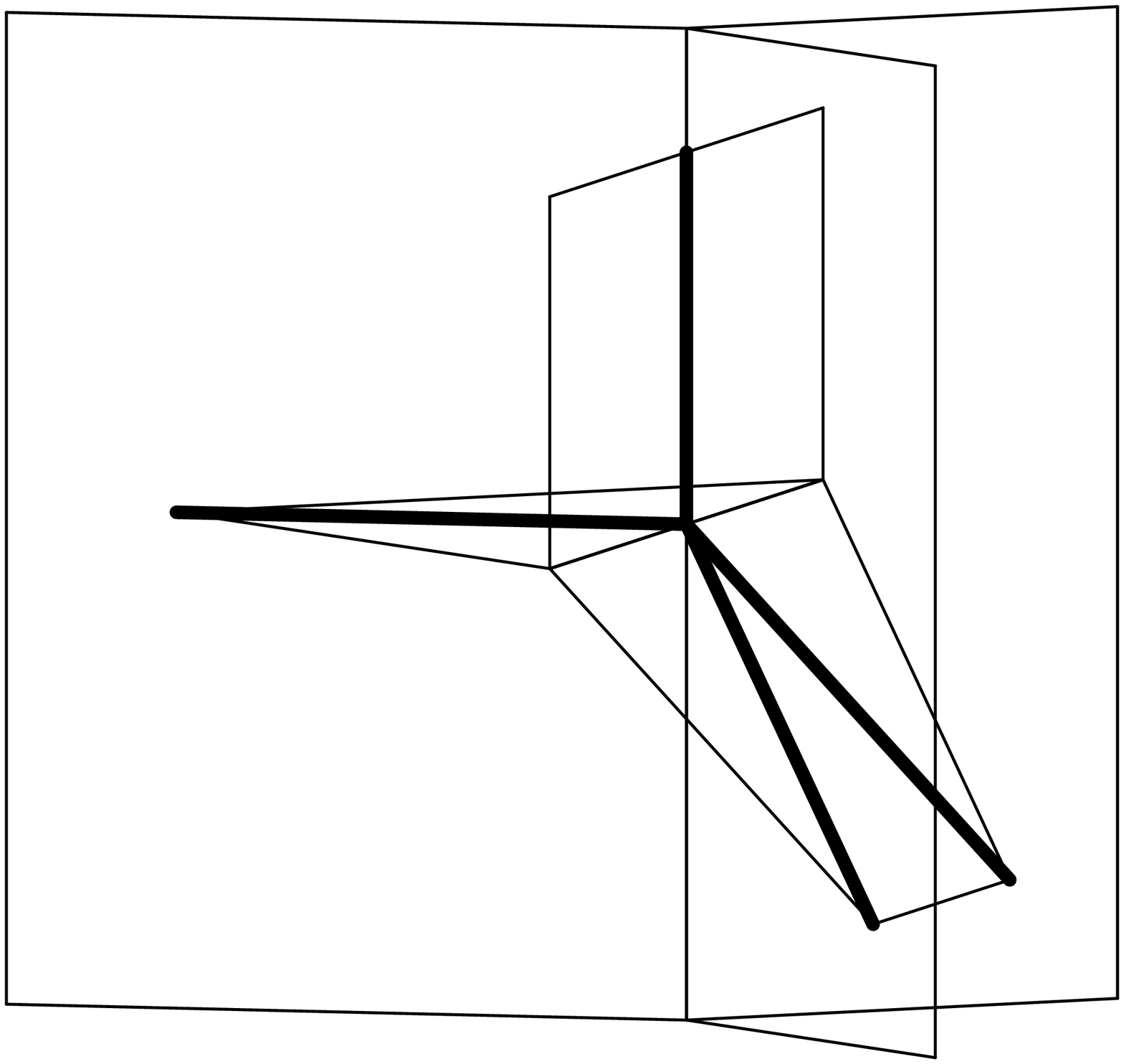}
\caption{A proper intersection which is not transversal.}\label{nontransv} 
\end{minipage}
\end{figure}

\begin{ex}
Figure \ref{3line} shows a {\em tropical line} in $\rr^3$ as the transversal intersection of two tropical planes (i.e., tropical hypersurfaces of degree 1).
\end{ex}

\begin{ex}
Figure \ref{nontransv} shows an intersection in $\rr^3$ which is proper, but not transversal. The surfaces are $X=V_{tr}(0x\oplus 0y\oplus 0)$ and $Y=V_{tr}(0xy\oplus 0z\oplus 0xyz)$. (Since the ``spines'' meet in a point, the intersection is not transversal.)
\end{ex}

\subsection{Intersection multiplicities}
Let $X_1,\dotsc,X_n\sub\rr^m$ be smooth tropical hypersurfaces such that the intersection $I=X_1\cap \dotsb\cap X_n$ is transversal. Let $U=X_1\cup \dotsb\cup X_n$ and denote by $\Subdiv_U$ the mixed subdivision associated to $U$. In \cite[Definition 4.3]{bert}, a general formula is given for the {\em intersection multiplicity} at each cell of $I$. For our purposes, two special cases suffice. If $P\in I^{(0)}$, let $P^\vee$ be the associated dual cell in $\Subdiv_U$.

\begin{dfn}\label{def:intmult1}
Suppose $n=m$, so $I$ consists of finitely many points. The intersection multiplicity at a point $P\in I$ is defined by $m_P=\vol(P^\vee)$.
\end{dfn}

\begin{remark}
This generalizes the standard definition of intersection multiplicities of tropical plane curves.
\end{remark}

\begin{dfn}\label{def:intmult2}
Suppose $n=m-1$, so $I$ is one-dimensional. The intersection multiplicity at a vertex $P\in I$ is defined by $m_P=2\vol(P^\vee)$.
\end{dfn}

\begin{remark}\label{rem:half}
It follows from the definition of transversality that $P^\vee$ has a privileged representation of the form $P^\vee=\La_1+\dotsb +\La_{n-1}+\Delta$, where each $\La_i$ is a primitive lattice interval, and $\Delta$ is a primitive lattice triangle. It follows from this that $\vol(P^\vee)$ is always a positive multiple of $\frac12$.
\end{remark}

\subsection{Tropical versions of Bernstein's Theorem and Bezout's Theorem}
Given polytopes $\Delta_1,\dotsc,\Delta_m$ in $\rr^m$, we consider the map $\ga:(\rr_{\geq 0})^m\pil \rr$ defined by $(\la_1,\dotsc,\la_m)\mapsto \vol(\la_1\Delta_1+\dotsb +\la_m\Delta_m)$. One can show that $\ga$ is given by a homogeneous polynomial in $\la_1,\dotsc,\la_m$ of degree $m$. We define the {\em mixed volume} of $\Delta_1,\dotsc,\Delta_m$ to be the coefficient of $\la_1\la_2\dotsm\la_m$ in the polynomial expression for $\ga$. The following tropical version of Bernstein's Theorem is proved in \cite[Corollary 4.7]{bert}:
\begin{theo}\label{theo:bern}
Suppose tropical hypersurfaces $X_1,\dotsc,X_m\sub \rr^m$ with Newton polytopes $\Delta_1,\dotsc,\Delta_m$ intersect in finitely many points. Then the total number of intersection points counted with multiplicities is equal to the mixed volume of $\Delta_1,\dotsc,\Delta_m$.
\end{theo}

As a special case of this we get a tropical version of Bezout's Theorem:
\begin{cor}\label{cor:bez}
Suppose the tropical hypersurfaces $X_1,\dotsc,X_m\sub \rr^m$ have degrees $d_1,\dotsc,d_m$, and intersect in finitely many points. Then the number of intersection points counting multiplicities is $d_1\dotsm d_m$.
\end{cor}
\begin{proof}
By Theorem \ref{theo:bern}, the number of intersection points, counting multiplicities, is the coefficient of $\la_1\la_2\dotsm\la_m$ in 
$$\vol(\la_1\Ga_{d_1}^m+\dotsb +\la_m\Ga_{d_m}^m)=\vol(\Ga_{\la_1d_1+\dotsb +\la_md_m}^m)=\frac1{m!}(\la_1d_1+\dotsb +\la_md_m)^m.$$
By the multinomial theorem, the wanted coefficient is $d_1\dotsm d_m$, as claimed.
\end{proof}

\section{Tropical complete intersection curves}
A {\em tropical complete intersection curve} $\cc$ is a transversal intersection of $n$ smooth tropical hypersurfaces $X_1,\dotsc,X_n\sub \rr^{n+1}$, for some $n\geq 2$. It is a one-dimensional polyhedral complex, some of whose edges are unbounded. We say that $\cc$ is {\em smooth} if the intersection multiplicity is 1 at each vertex (cf. Definition \ref{def:intmult2}).

Recall that any cell $C$ of $\cc$ is also a cell of the tropical hypersurface $U=X_1\cup\dotsb\cup X_n$. In particular, the notation $C^\vee$ always refers to the cell of $\Subdiv_U$ dual to $C\sub U$.

\begin{lem}\label{lem:3val}
Each vertex of $\,\cc$ has valence 3.
\end{lem}
\begin{proof}
If $P$ is a vertex of $\cc$, then by Remark \ref{rem:half}, $P^\vee$ has a privileged representation $P^\vee=\La_1+\dotsb +\La_{n-1}+\Delta$, where each $\La_i$ is a primitive interval, and $\Delta$ is a primitive lattice triangle. If $E$ is any edge of $\cc$ adjacent to $P$, then $E^\vee$ must be a mixed cell of $\Subdiv_U$ which is also a facet of $P^\vee$. This means that $E^\vee=\La_1+\dotsb +\La_{n-1}+\Delta'$, where $\Delta'$ is a side of $\Delta$. Hence there are exactly 3 such adjacent edges - one for each side of $\Delta$.
\end{proof}

Our first goal is to calculate the number of vertices of $\cc$. Before stating the general formula, let us discuss the easiest case as a warm up example:

\subsection{Example: Complete intersections in $\rr^3$}
Let $\cc=X\cap Y\sub \rr^3$ be a tropical complete intersection curve, where $X=V_{tr}(f)$ and $Y=V_{tr}(g)$ are smooth tropical surfaces of degrees $d$ and $e$ respectively.

\begin{theo}\label{theo:r3}
The number of vertices of $\cc$, counting multiplicities, is $de(d+e)$. 
\end{theo}
\begin{proof}
The idea is to look at all the vertices of the union $X\cup Y$, and their dual polytopes in the subdivision corresponding to $X\cup Y$. 
Since the intersection of $X$ and $Y$ is transversal, we can write the set of vertices of $X\cup Y$ as a disjoint union,
\begin{equation}\label{key3}
  (X\cup Y)^{(0)}= X^{(0)}\sqcup Y^{(0)}\sqcup (X\cap Y)^{(0)}.
\end{equation}
Now, any element $P\in (X\cup Y)^{(0)}$ corresponds to a maximal cell $P^\vee$ in $\Subdiv(f\odot g)$. The privileged representation of $P^\vee$ is of one of the following forms:
\begin{itemize}
\item $P^\vee=(\text{3-cell of $\Subdiv(f)$})+ (\text{0-cell of $\Subdiv(g)$})\quad\Longrightarrow\quad P\in X^{(0)}$. 
\item $P^\vee=(\text{0-cell of $\Subdiv(f)$})+ (\text{3-cell of $\Subdiv(g)$})\quad\Longrightarrow\quad P\in Y^{(0)}$. 
\item $P^\vee=(\text{2-cell of $\Subdiv(f)$})+(\text{1-cell of $\Subdiv(g)$})$ or\\ $P^\vee=(\text{1-cell of $\Subdiv(f)$})+ (\text{2-cell of $\Subdiv(g)$})\quad\Longrightarrow\quad P\in (X\cap Y)^{(0)}$.

\end{itemize}

Hence, dualizing \eqref{key3} and taking volumes, we get the relation
\begin{equation}\label{sum3}
  \sum_{P\in (X\cup Y)^{(0)}}\nnsp \vol(P^\vee)\; =\; \sum_{P\in X^{(0)}} \nnsp\vol(P^\vee)\;\;+\sum_{P\in Y^{(0)}} \nnsp\vol(P^\vee)\;\;+\sum_{P\in (X\cap Y)^{(0)}} \nnsp\vol(P^\vee).
\end{equation}

Now, if $P\in (X\cap Y)^{(0)}$, the volume of $P^\vee$ is $\frac12 m_P$ (by definition of intersection multiplicity).
Hence, \eqref{sum3} gives
\begin{equation*}
 \vol(\Delta_{f\odot g})=\vol(\Delta_f)+\vol(\Delta_g)+\sum_{P\in (X\cap Y)^{(0)}}\frac12 m_P.
\end{equation*}
Since $\Delta_f=\Ga_d^3$, $\Delta_g=\Ga_e^3$, and $\Delta_{f\odot g}=\Ga_d^3 + \Ga_e^3 =\Ga_{d+e}^3$, we find that
\begin{equation*}
\sum_{P\in (X\cap Y)^{(0)}} m_P= 2\Bigl[ \frac{(d+e)^3}{6}-\frac{d^3}{6}-\frac{e^3}{6}\Bigr] =de(d+e).
\end{equation*}
\end{proof}

\subsection{The number of vertices in the general case}
In this section we prove the following generalization of Theorem \ref{theo:r3}:

\begin{theo}\label{nodetheo}
Let $\cc=X_1\cap\dotsb \cap X_n$ be a tropical complete intersection curve in $\rr^{n+1}$, where $X_1,\dotsc,X_n$ are smooth of degrees $d_1,\dotsc, d_n$. The number of vertices of $\cc$, counting multiplicities, is
\begin{equation*}
\sum_{P\in \cc^{(0)}}m_P =d_1d_2\dotsb d_n(d_1+d_2+\dotsb +d_n).
\end{equation*}
\end{theo}

To prove Theorem \ref{nodetheo}, we will use the same setup as in the previous section. Note that in the proof of the case $n=3$, the relation \eqref{key3} is the key giving us control over $(X\cap Y)^{(0)}$. So as an auxiliary lemma, we first state and prove a generalization of this. 

To simplify the writing, we introduce the following notation: Let $[n]=\{1,2,\dotsc,n\}$. For any nonempty subset $J=\{j_1,\dotsc,j_k\}\sub[n]$, we put
\begin{equation}\label{eq:UI}
\begin{split}
U_J& :=X_{j_i}\cup \dotsb\cup X_{j_k},\\
I_J& :=X_{j_i}\cap \dotsb\cap X_{j_k}.
\end{split}\end{equation}
In the special case $J=[n]$, we simply write $U$ and $I$, i.e. $U:=U_{[n]}$ and $I=\cc=I_{[n]}$.

By the assumption of transversality, we have $I_{\!J}^{\,(0)}\cap I_{\!K}^{\,(0)}=\emptyset$ whenever $J, K\sub [n]$ are distinct nonempty subsets. Thus we can split the $0$-cells of $U=X_1\cup \dotsb\cup X_n$ into a disjoint union:
\begin{equation*}
U^{(0)}=\bigsqcup_{J\sub [n]} I_{\!J}^{\,(0)}.
\end{equation*}
Similarly, for any nonempty subset $J\sub[n]$, we get
\begin{equation}\label{nodeUnionSub}
\ujn=\bigsqcup_{J'\sub J} \ijmn.
\end{equation}

\begin{lem}\label{keylem}
For a transversal intersection of tropical hypersurfaces $X_1,\dotsc, X_n$, we have:
\begin{equation}\label{keylemexpr}
I^{(0)}\;\sqcup\nsp \bigsqcup_{|J|=n-1}\nsp\nsp \ujn\;\sqcup\nsp \bigsqcup_{|J|=n-3}\nsp\nsp \ujn\;\sqcup\dotsb\\
=U^{(0)}\;\sqcup \nsp \bigsqcup_{|J|=n-2}\nsp\nsp \ujn\;\sqcup\nsp \bigsqcup_{\;|J|=n-4}\nsp\nsp \ujn\;\sqcup\dotsb.\nsp\nsp
\end{equation}
\end{lem}

\begin{proof}
By applying \eqref{nodeUnionSub} to every set $\ujn$ in \eqref{keylemexpr}, we see that the following expression is equivalent to \eqref{keylemexpr}:
\begin{equation}\label{keylemexpr2}
I^{(0)}\;\sqcup \nsp\bigsqcup_{\substack{|J|=n-1\\J'\sub J}}\nsp\nsp \ijmn\;\sqcup\nsp \bigsqcup_{\substack{|J|=n-3\\J'\sub J}}\nsp\nsp \ijmn\;\sqcup\dotsb
=\bigsqcup_{J'\sub [n]}\nsp \ijmn\;\sqcup\nsp \bigsqcup_{\substack{|J|=n-2\\J'\sub J}}\nsp\nsp \ijmn\;\sqcup\nsp \bigsqcup_{\substack{|J|=n-4\\J'\sub J}}\nsp\nsp \ijmn\;\sqcup\dotsb.\nsp\nsp\nsp\nsp
\end{equation}

We claim that for each fixed subset $J'\sub [n]$, the set $\ijmn$ appears equally many times on each side of \eqref{keylemexpr2}. By inspection, this is true for $J'=[n]$. Assume now $|J'|=k< n$. Then for any integer $s$ with $k\leq s\leq n$, there are exactly $\binom{n-k}{s-k}$ sets $J\sub [n]$ containing $J'$ such that $|J|=s$. Hence, the number of times $\ijmn$ appears on the left side of \eqref{keylemexpr2} is $\binom{n-k}{n-1-k}+\binom{n-k}{n-3-k}+\dotsb = \binom{n-k}{1}+\binom{n-k}{3}+\dotsb =2^{n-k-1}$,
while the number of appearances on the right side is $\binom{n-k}{n-k}+\binom{n-k}{n-2-k}+\dotsb = \binom{n-k}{0}+\binom{n-k}{2}+\dotsb =2^{n-k-1}$. This proves the claim, and the lemma follows.
\end{proof}

\begin{proof}[Proof of Theorem \ref{nodetheo}]
Suppose $\cc$ and $X_1,\dotsc, X_n$ are as in the statement of the theorem. We assume that for each $i$, $X_i$ has degree $d_i$, so the associated Newton polytope is the simplex $\Ga_{d_i}^{n+1}$.
Let $U$ denote the union $X_1\cup\dotsb\cup X_n$, and $\Subdiv_U$ the associated subdivision of 
$\Ga_{d_1+\dotsb +d_n}^{n+1}$. 

For each nonempty $J=\{j_1,\dotsc,j_k\}\sub [n]$, let $U_J$ and $I_J$ be as in \eqref{eq:UI}. In particular, $U_J$ is a tropical hypersurface (set-theoretically contained in $U$) with an associated subdivision $\Subdiv_{U_J}$ of the simplex $\Delta_J:=\Ga^{n+1}_{d_{j_1}+\dotsb +d_{j_k}}$. 

Each vertex of $U_J$ is also a vertex of $U$, and therefore corresponds to a maximal cell of $\Subdiv_U$. Let $S_J$ be the set of maximal cells of $\Subdiv_U$ corresponding to the vertices of $U_J$.
By transversality, the elements of $S_J$ are simply translations of the maximal cells of $\Subdiv_{U_J}$. Hence the total volume of the cells of $S_J$, denoted $\vol(S_J)$, is 
\begin{equation*}
\vol(S_J)=\sum_{P\in \ujn}\vol(P^\vee)=\vol(\Delta_J)=\frac1{(n+1)!}(d_{j_1}+\dotsb +d_{j_k})^{n+1}.
\end{equation*}

Now we turn to Lemma \ref{keylem}. Dualizing \eqref{keylemexpr}, we find that
\begin{equation}
\sum_{P\in I^{(0)}}\nsp\nsp \vol(P^\vee)\; + \sum_{|J|=n-1}\nsp\nsp \vol(S_J)\; + \dotsb\;
=\; \vol(S) \:+ \sum_{|J|=n-2}\nsp\nsp \vol(S_J)\;+\dotsb
\end{equation}
By the definition of intersection multiplicity, the dual $P^\vee\in \Subdiv_U$ of a vertex $P\in I^{(0)}$ has volume $\frac12 m_P$. It follows that
\begin{equation*}
  \sum_{P\in I^{(0)}}\frac12 m_P =\frac1{(n+1)!}\sum_{\{j_i,\dotsc,j_k\}\sub[n]} (-1)^{n-k}(d_{j_1}+\dotsb + d_{j_k})^{n+1},
\end{equation*}
which after some elementary manipulation reduces to $$\sum_{P\in I^{(0)}} m_P=d_1d_2\dotsm d_n(d_1+d_2+\dotsb +d_n).$$
\end{proof}

\subsection{The genus of tropical complete intersection curves}

\begin{dfn}
The {\em genus} $g=g(\cc)$ of a tropical complete intersection curve $\cc$ is the first Betti number of $\cc$, i.e., the number of independent cycles of $\cc$.
\end{dfn}

\begin{lem}\label{lem:genus}
For a connected tropical complete intersection curve $\cc$, we have
\begin{equation*}
  2g(\cc)-2=v-x,
\end{equation*}
where $v$ is the number of vertices, and $x$ the numbers unbounded edges of $\cc$.
\end{lem}
For the proof, recall that a graph is called {\em 3-valent} if every vertex has 3 adjacent edges. Furthermore, we apply the following terminology: A one-dimensional polyhedral complex in $\rr^m$ with unbounded edges is regarded as a graph, where the 1-valent vertices have been removed. For example, a tropical line in $\rr^3$ is considered a 3-valent graph with 2 vertices and 5 edges.
\begin{proof}
By Lemma \ref{lem:3val}, $\cc$ is 3-valent. Since $\cc$ is connected, it has a spanning tree $T$, such that $\cc\smallsetminus T$ consists of $g$ edges. While $T$ is not 3-valent, we can construct a 3-valent tree $T'$ from $T$ by adding unbounded edges wherever necessary. Clearly, we must add exactly $2g$ such edges. Thus if $\cc$ has $v$ vertices and $e$ edges, $T'$ has $v$ vertices and $e+g$ edges. Since $T'$ is 3-valent, it is easy to see (for example by induction) that the number of edges is one more that twice the number of vertices, i.e., 
\begin{equation}\label{eq:edges}
e+g-1=2v.
\end{equation}
On the other hand, since $\cc$ is 3-valent, we must have $e=\frac12 (3v+x)$. Combining this with \eqref{eq:edges} gives the wanted result. 
\end{proof}

\begin{lem}\label{lem:unbounded}
Let $\cc$ be the transversal intersection of $X_1,\dotsc,X_n\sub\rr^{n+1}$, where each $X_i=V_{tr}(f_i)$ is a smooth tropical hypersurface of degree $d_i$. If $\cc$ is smooth, the number of unbounded edges of $\cc$ is $x=(n+2)d_1\dotsm d_n$.
\end{lem}
\begin{proof}
Let $U=X_1\cup \dotsb\cup X_n$, and let $\Subdiv_U$ be the associated subdivision of the simplex $\Ga:=\Ga_{d_1+\dotsb +d_n}^{n+1}$. The unbounded edges of $\cc$ are then in one-one correspondence with the mixed $n$-cells of $\Subdiv_U$ contained in the boundary of $\Ga$. To prove the lemma, it therefore suffices to show that there are exactly $d_1\dotsm d_n$ mixed $n$-cells in each of the $n+2$ facets of $\Ga$. By symmetry it is enough to consider the facet $\Ga'$ with $e_1=(1,0,\dotsc,0)$ as an inner normal vector. In the following we identify $\rr^n$ with the hyperplane in $\rr^{n+1}$ orthogonal to $e_1$.

For each $i=1,\dotsc, n$ let $\mathcal{S}_i$ be the subdivision induced by $\Subdiv(f_i)$ on the facet of $\Ga_{d_i}^{n+1}$ with $e_1$ as an inner normal vector. We can then regard $\mathcal{S}_i$ as the subdivision associated to the tropical hypersurface $X'_i:=V_{tr}(f'_i)\sub \rr^n$, where $f'_i$ is the tropical polynomial obtained from $f_i$ by removing all terms containing $x_1$. Furthermore, $X'_i$ is homeomorphic to the intersection $X_i\cap H$, where $H$ is any (classical) hyperplane with equation $x_1=k$ and $k<<0$. Note that $\deg X'_i=\deg X_i=d_i$.

Let $\mathcal{S}$ be the subdivision of $\Ga'$ induced by $\Subdiv_U$. As above, we regard $\mathcal{S}$ as the subdivision associated to the union $X'_1\cup \dotsb \cup X'_n\sub\rr^n$. Thus, the (finitely many) points in the intersection $I:=X'_1\cap \dotsb \cap X'_n$ are precisely the duals of the mixed $n$-cells of $\mathcal{S}$. 
We know from Theorem \ref{cor:bez} that the number of points in $I$ is $d_1\dotsm d_n$ when counting with intersection multiplicities; in other words (by Definition \ref{def:intmult1}) we have $\sum_{Q\in I} \vol_n(Q^\vee)=d_1\dotsm d_n$. 

All that remains is to show that if $Q\in I$, then $\vol_n(Q^\vee)=1$. This is where smoothness of $\cc$ comes in: Regarding $Q^\vee$ as an $n$-cell in $\Subdiv_U$, let $P$ be the vertex of $\cc$ such that $Q^\vee$ is a facet of $P^\vee\in\Subdiv_U$. Writing (as in Remark \ref{rem:half}) $P^\vee=\La_1+\dotsb +\La_{n-1}+\Delta$, where the $\La_i$'s are primitive intervals and $\Delta$ a primitive triangle, we must have $Q^\vee=\La_1+\dotsb +\La_{n-1}+\Delta'$, where $\Delta'$ is a side in $\Delta$. Since $\vol_{n+1}(P^\vee)=\frac12$ (by smoothness), it follows from this that $\vol_n(Q^\vee)=1$.
\end{proof}

\begin{theo}\label{theo:genus}
Let $\cc$ be the transversal intersection of $n$ smooth tropical hypersurfaces in $\rr^{n+1}$ of degrees $d_1,\dotsc, d_n$. If $\cc$ is smooth and connected, the genus $g$ of $\cc$ is given by
\begin{equation}\label{eq:g}
2g-2= d_1\dotsm d_n(d_1+\dotsb + d_n-(n+2)).
\end{equation}
\end{theo}
\begin{proof}
Since $\cc$ is smooth, it has exactly $v=d_1d_2\dotsm d_n(d_1+d_2+\dotsb +d_n)$ vertices (by Theorem \ref{nodetheo}) and $x=(n+2)d_1\dotsm d_n$ unbounded edges (by Lemma \ref{lem:unbounded}). Combined with Lemma \ref{lem:genus}, this proves the theorem.
\end{proof}

\begin{remark}
In complex projective space it is well known that any complete intersection curve is connected. This follows from standard cohomological arguments (see also \cite[Section 3.4.6]{Harts_connect} for a direct geometric argument due to Serre). In the tropical setting, it is known that any transversal intersection of tropical hyperplanes is a {\em tropical variety}, i.e., the tropicalization of an algebraic variety defined over the field of Puiseux series (\cite[Section 3, and Lemma 1.2 for the relation to Puiseux series]{Computing}). Furthermore, if a tropical variety is the tropicalization of an irreducible variety, then it is connected (\cite[Theorem 2.2.7]{EKL}). This suggests that - at least in the general case - a tropical complete intersection curve is connected. However, to the author's knowledge, this has not been proved.
\end{remark}

\begin{remark}
The formula \eqref{eq:g} coincides with the genus formula for a smooth complete intersection in $\pp^{n+1}_{\mathbb{C}}$ of $n$ hypersurfaces of degrees $d_1,\dotsc, d_n$.
\end{remark}

\section{Example: Tropical elliptic curves in $\rr^3$}
\begin{figure}[tb]
\centering
\includegraphics[scale=0.3]{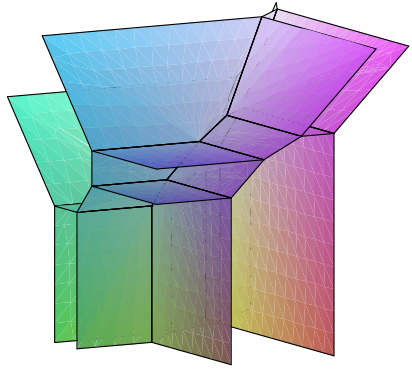}\label{Q}\caption{A tropical quadric surface in $\rr^3$.}
\end{figure}

By a {\em tropical quadric surface} in $\rr^3$, we mean a smooth tropical hypersurface of degree 2.  
In this section we take a closer look at intersections of {\em tropical quadric surfaces} in $\rr^3$, i.e., smooth tropical hypersurfaces in $\rr^3$ of degree 2. Figure \ref{Q} shows a typical tropical quadric surface. 

Let $\cc$ be a smooth, connected complete intersection curve of two tropical quadric surfaces in $\rr^3$. We call $\cc$ a {\em tropical elliptic curve}. The name is justified by Theorem \ref{theo:genus}, which tells us the the genus $g$ of $\cc$ satisfies $2g-2= 2\cdot 2\cdot (2+2-4)$, that is, $g=1$. In particular, $\cc$ contains a unique cycle.

Since $\cc$ is smooth, it has exactly $2\cdot 2\cdot (2+2)=16$ vertices, by Theorem \ref{nodetheo}. We divide these into two categories: Those on the cycle (called {\em internal vertices}), and the rest ({\em external vertices}). Clearly, $\cc$ has at least 3 internal vertices. But what is the maximum number of internal vertices? As the following example shows, all 16 vertices can be internal:

\begin{figure}[tb]
\noindent
\begin{minipage}[b]{.48\linewidth}
\centering
\includegraphics[scale=0.2]{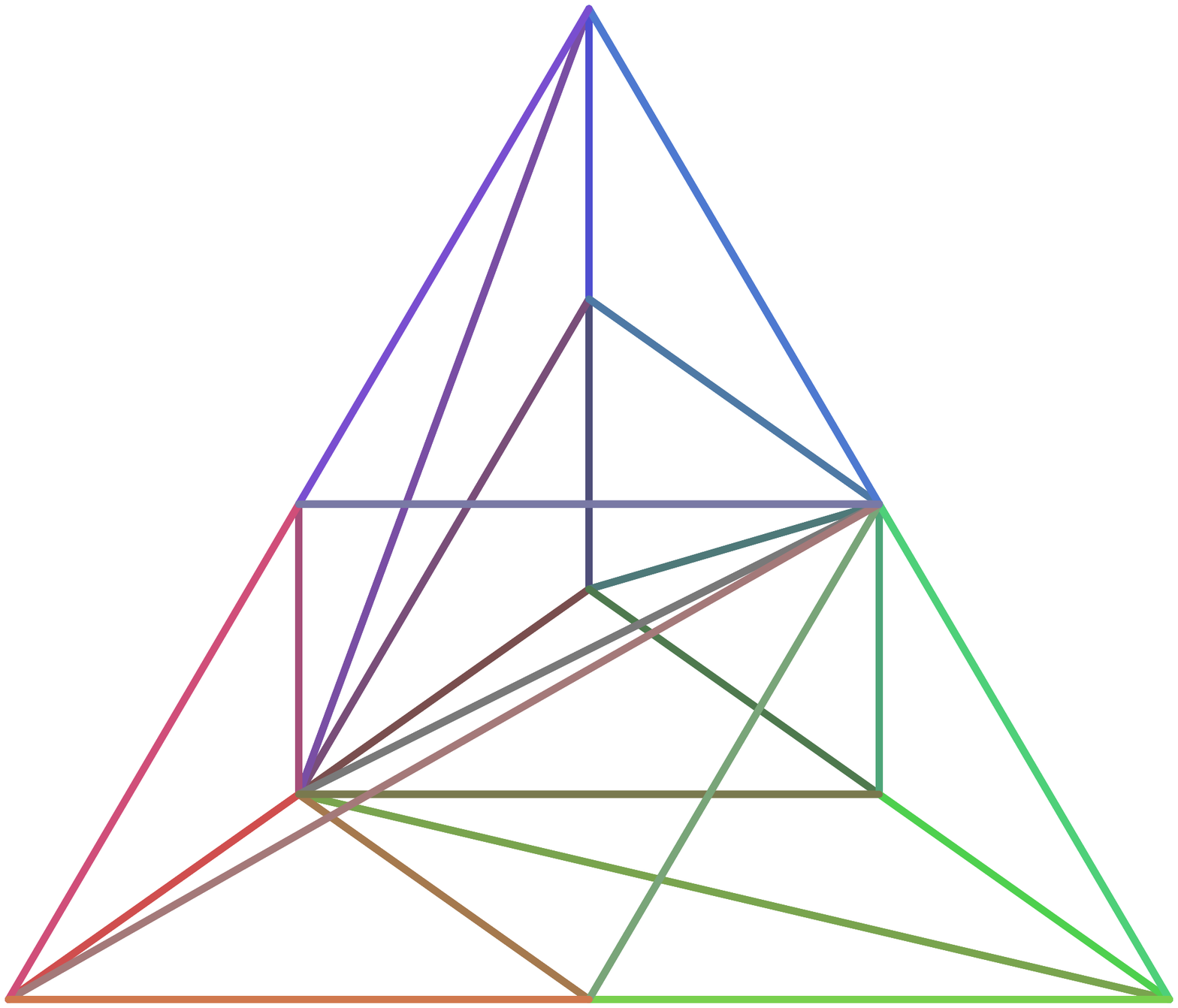}\caption{The subdivision $\Subdiv(f)$.}\label{Q1sub}
\end{minipage}\hfill
\begin{minipage}[b]{.48\linewidth}
  \centering
\includegraphics[scale=.2]{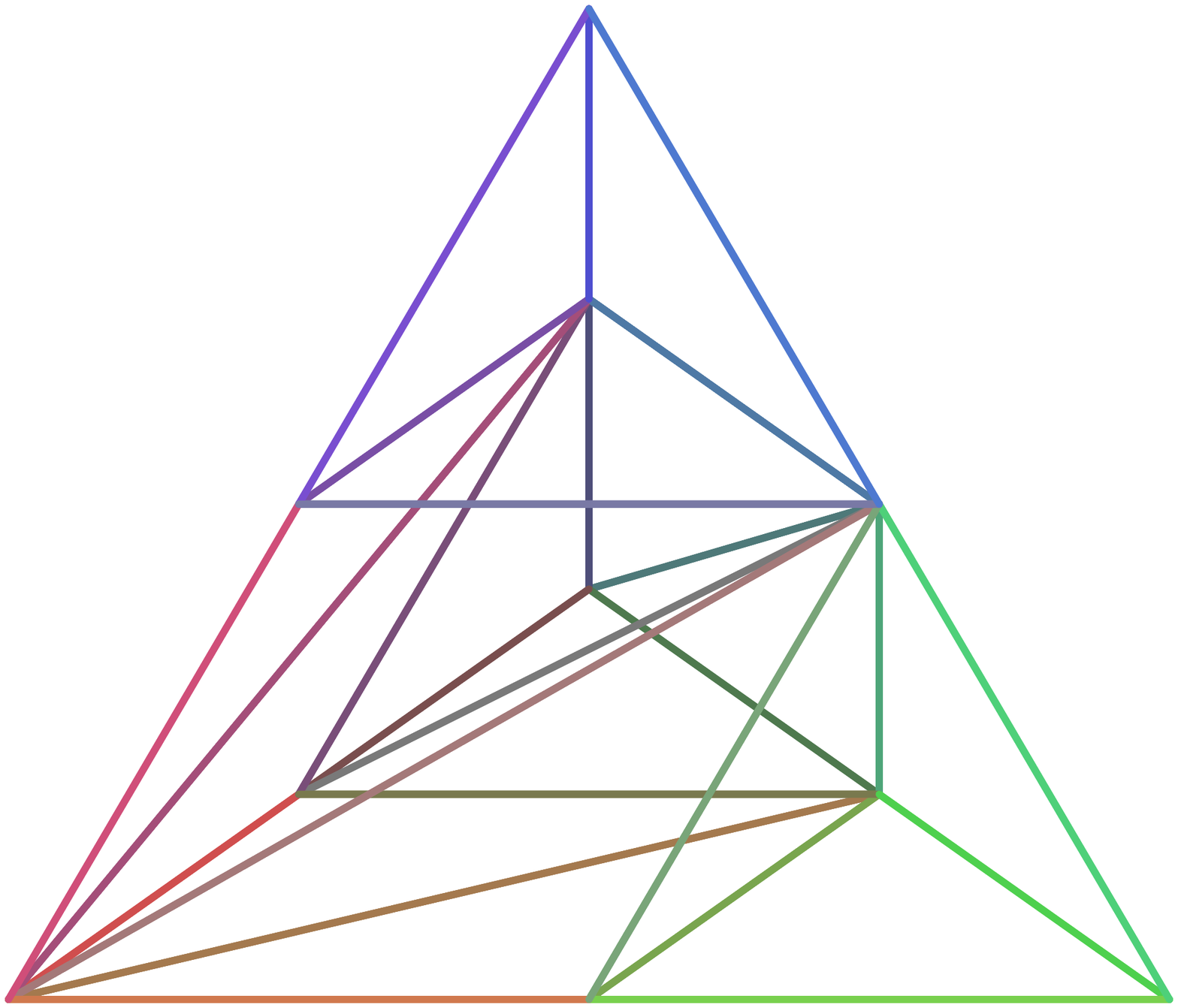}\caption{The subdivision $\Subdiv(g)$.}\label{Q2sub}
\end{minipage}
\end{figure}

\begin{figure}[tb]
\noindent
\begin{minipage}[b]{.48\linewidth}
\centering
\includegraphics[scale=0.3]{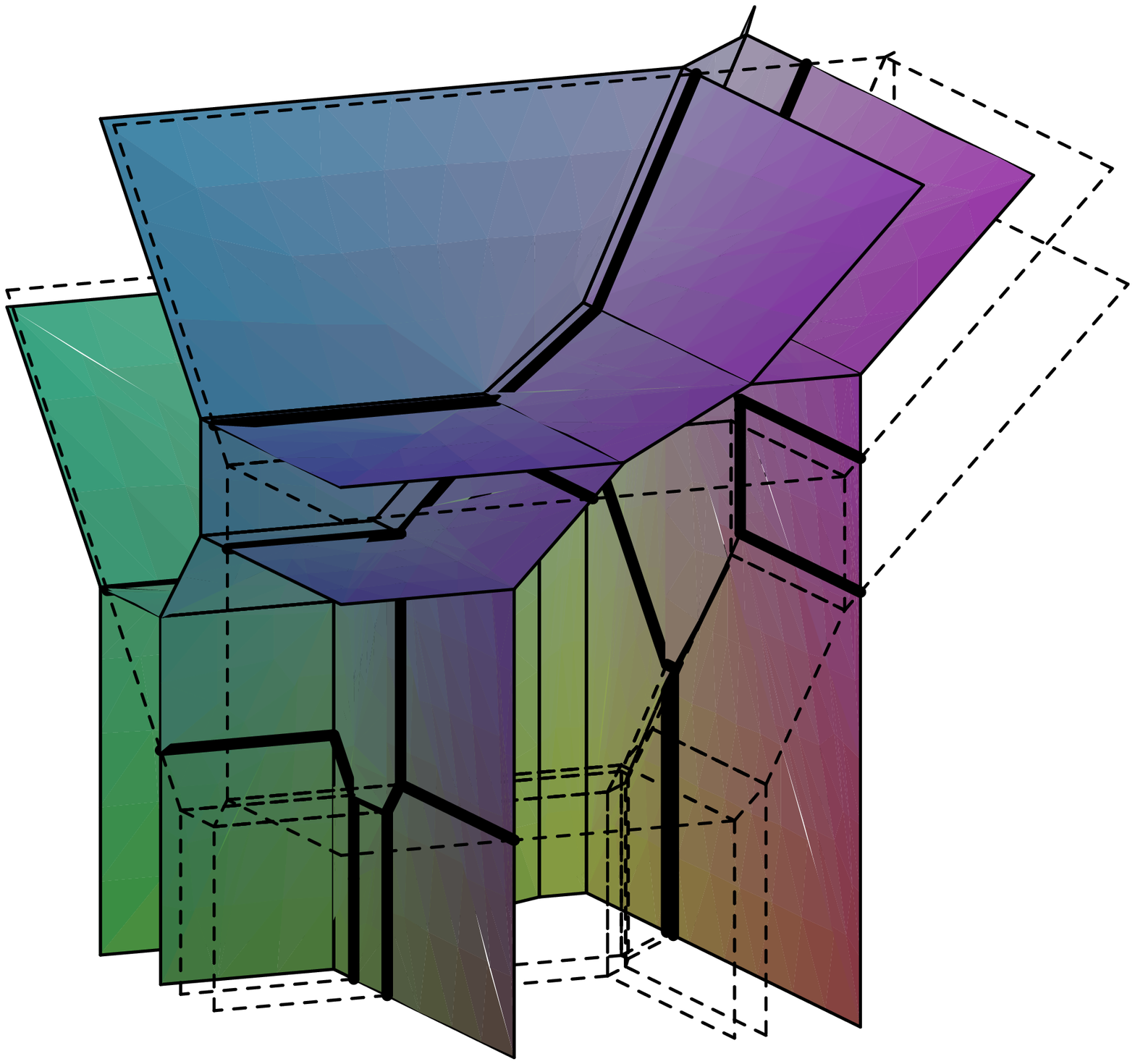}\caption{The intersection $\cc=Q_1\cap Q_2$.}\label{recint}
\end{minipage}\hfill
\begin{minipage}[b]{.48\linewidth}
  \centering
\includegraphics[scale=.3]{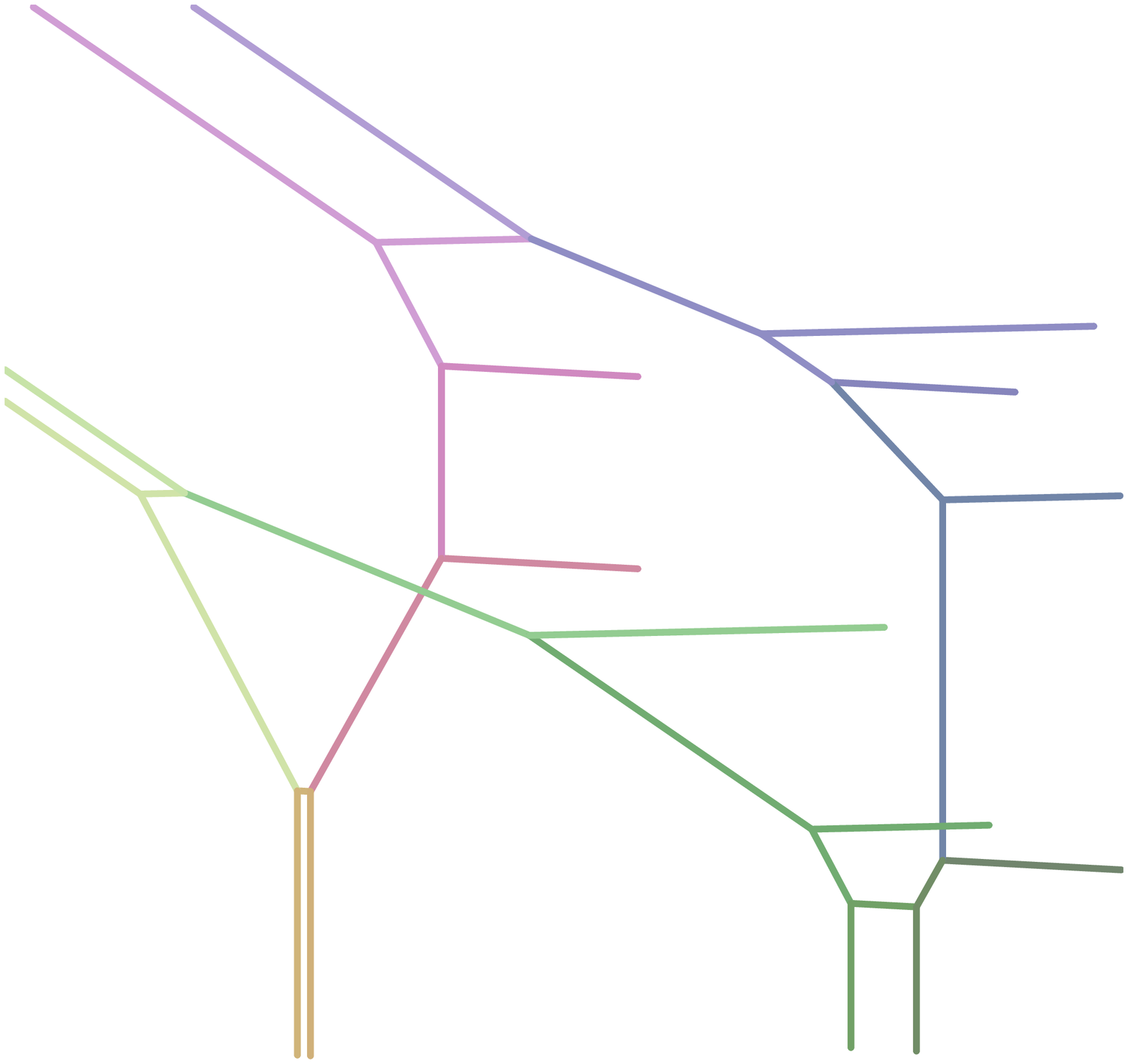}\caption{$\cc$ has 16 internal vertices.}\label{rec}
\end{minipage}
\end{figure}

\begin{ex}
Let $Q_1=V_{tr}(f)$ and $Q_2=V_{tr}(g)$, where
\begin{multline*}
f(x,y,z)=(-6)\oplus 13x\oplus (-3)y\oplus (-4)z\oplus 10x^2\oplus 2xy
\oplus 4xz\oplus (-9)y^2\oplus 5yz\oplus (-9)z^2,
\end{multline*}
and
\begin{multline*}
g(x,y,z)=(-15)\oplus (-10)x\oplus (-4)y\oplus 2z\oplus (-7)x^2\oplus (-2)xy\\
\oplus 0xz\oplus 2y^2\oplus 15yz\oplus (-1)z^2.
\end{multline*}
Figures \ref{Q1sub} and \ref{Q2sub} show the subdivisions of $\Ga_2^3$ induced by $f$ and $g$ respectively.

The intersection curve $\cc=Q_1\cap Q_2$ has genus 1 and 16 internal vertices. Figure \ref{recint} shows the the two quadrics intersecting. In Figure \ref{rec} we see the intersection curve alone from a different angle, clearly showing the cycle with all its 16 vertices.
\end{ex}

\begin{remark}
A computer search shows that for every integer $m$, with $3\leq m\leq 16$, there exist two tropical quadric surfaces in $\rr^3$ intersecting transversally in a tropical elliptic curve with $m$ internal vertices. 
\end{remark}
\vspace{.2cm}
\noindent
{\it Acknowledgements}.
I would like to thank my supervisor Kristian Ranestad for his constant support. I am also grateful to Bernd Sturmfels for valuable discussions about the material of this paper, and for posing the problem of how many internal vertices a tropical elliptic curve can have. 

\bibliographystyle{plain}
\bibliography{bib_complete}
\end{document}